\documentclass{amsart}

\newtheorem{theorem}{Theorem}[section]
\newtheorem{lemma}{Lemma}[section]
\newtheorem{proposition}{Proposition}[section]

\theoremstyle{definition}
\newtheorem{definition}{Definition}[section]
\newtheorem{example}{example}[section]

\theoremstyle{remark}
\newtheorem{remark}{Remark}[section]

\numberwithin{equation}{section}

\usepackage{amssymb}

\usepackage[shortlabels]{enumitem}
\setlist[enumerate,1]{label={\upshape(\roman*)}}

\renewcommand\leftmark{L. LAMGOUNI}
\renewcommand\rightmark{SPHERICAL-VECTORS}

\usepackage{textcmds}

\usepackage{hyperref}

\usepackage[capitalise]{cleveref}

\usepackage{graphicx}
\graphicspath{ {./Figs/} }
\usepackage{caption}
\usepackage{float}

\begin{document}

\title{Spherical-vectors and geometric interpretation of unit quaternions}

\author{Lahcen Lamgouni\textsuperscript{*}}

\curraddr{Minister of National Education Preschool and Sports: Errachidia, Morocco}
\email{lahcen.lamgouni@gmail.com (\href{https://orcid.org/0000-0002-3927-0353}{https://orcid.org/0000-0002-3927-0353}).}


\subjclass[2020]{Primary 51N30, 51A25; Secondary, 51N99}



\keywords{quaternions, quaternion multiplication, geometric interpretation, argument, polar form, exponential form, visualizing quaternions, non-commutativity, spherical form, spherical-vectors.}

\maketitle

\begin{center}
\noindent
\textsuperscript{*}\href{https://orcid.org/0000-0002-3927-0353}{https://orcid.org/0000-0002-3927-0353}\\
Email: \texttt{lahcen.lamgouni@gmail.com}
\end{center}

\begin{abstract}
In this article, we introduce and study the concept of \textit{spherical-vectors}, which can be perceived as a natural extension of the arguments of complex numbers in the context of quaternions. We initially establish foundational properties of these distinct vectors, followed by a demonstration, through the transfer of structure, that spherical-vectors constitute a non-abelian additive group, isomorphic to the group of unit quaternions. This identification facilitates the presentation of a novel polar form of quaternions, highlighting its algebraic properties, as well as the algebraic properties of the exponential writing. Furthermore, it enables the depiction of unit quaternions on the unit sphere of $\mathbb{R}^3$, allowing for a geometric interpretation of their multiplication.
\end{abstract}

\section{Introduction}\label{sec:sec1}

The algebraic properties of the additive group of oriented angles between two vectors in the Euclidean plane have allowed for the definition of the argument of complex numbers and the utilization of the algebraic properties of their exponential form to their full advantage. However, if we move to the three-dimensional Euclidean space, most of the preceding results collapse. Thus, There is no way to define the measure of an oriented angle between two vectors (in the case where we consider this measure to be a real number), and the arguments of quaternions are non-oriented and unstructured angles that no longer follow the Chasles relation. As T.Y. Lam wrote in his book \cite{Lam} \qq{\textit{The failure of the equality $e^{p}\cdot e^{q}= e^{p+q}$ on $\mathbb{H}$ is a serious drawback, and may have been the principal roadblock to the development of a really useful theory of functions on the quaternions}}. For more details on the argument and the polar form of quaternions, please refer to \cite{Kuipers,Lam,Quaternions}.

\begin{figure}[H]%
\centering
\includegraphics[width=1\textwidth]{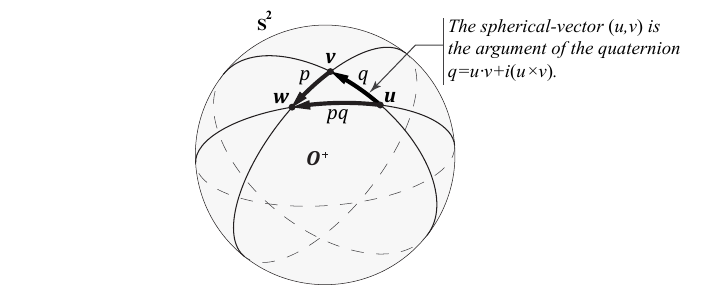}
\caption{The multiplication $pq$ of two unit quaternions $p=v\cdot w+i(v\times w)$ and $q=u\cdot v+i(u\times v)$ can be geometrically interpreted through the elementary Chasles relation $(u,v)+(v,w)=(u,w)$.}\label{fig:fig1}
\end{figure}

The present work aims to generalize, in a way naturally compatible with complex numbers, the notions of, argument, polar form and exponential writing to quaternions. Thus, we obtain a manageable tool allowing to establish algebraically and geometrically all the properties related to the multiplication of quaternions in their exponential forms, namely the fundamental equation\footnote{We will see later that $\alpha$ and $\beta$ represent the new concept of argument for quaternions, and that their addition is non-commutative. Calculations will also demonstrate that, in the multiplication of quaternions, it is generally necessary to commute the addition of their new arguments.} $e^{i\alpha}e^{i\beta}=e^{i(\beta+\alpha)}$. To achieve this, we first need to define and study the notion of \textit{spherical-vectors}. This goes back to Hamilton's model for the group of unit quaternions (Hamilton \cite{Hamilton}, Book II, Chapter I, Section 9) and (\cite{Viro}, Subsection 6.3, page 22). In his model, Hamilton interpreted unit quaternions as vector-arcs on the unit sphere $S^2$ of the three euclidean space. In this paper, we reinterpret these vector-arcs from an alternative perspective to endow them with an additive algebraic structure. Consequently, we construe each unit quaternion $q$ as an ordered pair $(u,v)$ of unit vectors in $\mathbb{R}^3$, which we designate as \textit{spherical-vectors} and which serve as a novel representation of the argument of $q$. Following Hamilton's approach, the quaternion $q=x+yi+zj+tk$, with $(x,y,z,t)\in\mathbb{R}^4$, is expressed in the form $q=x+U$, where $U(y,z,t)$ is the vector component of $q$ depicted by the pure quaternion $U=yi+zj+tk$. From our standpoint, we represent $q$ in the slightly different way $q=x+iV$, which we call \textit{spherical form} of $q$, where $V$ is the vector $(t,-z,y)$ represented by the quaternion $V=tj-zk+y$. Throughout this article, we have adopted this novel vector convention to render the form of $q$ more coherent and in alignment with the algebraic and exponential forms of complex numbers, albeit with the potential to initially perplex the reader. Given that $x^2+\Vert V\Vert^2=1$, it is straightforward to identify a spherical-vector $(u,v)$, representing our novel argument of $q$, such that $q=u\cdot v+i(u\times v)$, with $u\cdot v=x$ and $u\times v=V$. Consequently, we introduce the exponential notation $q=e^{i(u,v)}$ with $\arg(q)=(u,v)$, and we demonstrate the fundamental algebraic property of the argument and exponential notation as follows. Given two unit quaternions $p$ and $q$ with arguments represented by spherical-vectors $\alpha$ and $\beta$, respectively, we show that it is possible to find three unit vectors $u$, $v$, and $w$ such that $\alpha=(v,w)$ and $\beta=(u,v)$. We then establish the result: \[pq=\left[v\cdot w+i(v\times w))\right]\left[u\cdot v+i(u\times v)\right]=u\cdot w+i(u\times w),\] which can be expressed as \[e^{i(v,w)}e^{i(u,v)}=e^{i\left[(u,v)+(v,w)\right]}=e^{i(u,w)}.\] In other words, $\arg(pq)=\arg(q)+\arg(p)$. In this case, arguments possess an additive algebraic structure on $S^2$ that is non-commutative (\Cref{fig:fig1}).
In this article, we propose an innovative approach compared to previous studies \cite{Goldman, Baek, Hanson}, which all address the issue of the geometric interpretation of quaternion multiplication based on rotations. The two main contributions of our article lie in the visualization of unit quaternions and the algebraic properties of quaternion multiplication in their exponential forms.

Our original method allows for the apprehension of unit quaternions and their multiplication through the use of spherical-vectors and their addition represented on $S^2$, without resorting to rotations. This approach highlights the algebraic properties of quaternion multiplication in their exponential forms, thus providing a new perspective on the understanding of these mathematical objects.

Our contribution is all the more significant as it extends the current knowledge of quaternions and their algebraic properties, thereby offering new perspectives for researchers and practitioners in this field.

We briefly outline the structure of the article. Section \ref{sec:sec2} is dedicated to the discussion of spherical-vectors. Initially, we present the definition and the primary properties associated with this concept. Subsequently, we offer a geometric representation of a spherical-vector on $S^2$. To establish a connection between spherical-vectors and unit quaternions, we introduce the spherical form of quaternions in Section \ref{sec:sec3}. In Section \ref{sec:sec4}, we construct a bijective map $\mu$ from the set $\mathcal{S}$ of spherical-vectors to the multiplicative group $\mathbb{H}_1$ of unit quaternions. This bijection enables us to transport the structure of $\mathbb{H}_1$ into $\mathcal{S}$. Consequently, we define the following additive internal composition law in $\mathcal{S}$: $\alpha + \beta = \mu^{- 1}\left(\mu (\beta) \mu(\alpha)\right)$. Equipped with this operation, $\mathcal{S}$ forms a non-commutative additive group, and $\mu$ serves as an anti-isomorphism of groups. Utilizing spherical-vectors, we propose in Section \ref{sec:sec5} a new reformulation of the definitions of argument and polar form of quaternions. We then demonstrate the algebraic properties of the argument and the exponential form. Lastly, Section \ref{sec:sec6} provides several examples of applications.

In the sequel the space $\mathbb{R}^3$ is endowed with its standard oriented Euclidean structure and its standard basis $(e_x,e_y,e_z)$. We designate respectively by \qq{$\ \cdot\ $}, \qq{$\ \times\ $} and \qq{$\ \Vert\cdot\Vert\ $}, the dot product, the cross product and the Euclidean norm.

\section{Spherical-Vectors}\label{sec:sec2}

In this section, we introduce the concept of a spherical-vector, detailing its components and support. We then proceed to illustrate its geometric representation on the unit sphere $S^2$.

\subsection{Definitions}\label{subsec:subsec2.1}

\begin{definition}\label{def:def2.1} [\textbf{Components of an ordered pair of non-zero vectors}]\\
Let $(u,v)$ be an ordered pair of non-zero vectors in $\mathbb{R}^3$. We call \textit{scalar component} and \textit{vector component} of $(u, v)$, the real number $\lambda = \frac{u \cdot v}{\Vert u \Vert \Vert v \Vert }$ and the vector $n = \frac{u \times v}{\Vert u \Vert \Vert v \Vert}$, respectively.
\end{definition}

\begin{definition}\label{def:def2.2}[\textbf{Spherical-vector}]
A \textit{spherical-vector} of $\mathbb{R}^3$ is an equivalence class of ordered pairs of non-zero vectors. Two ordered pairs represent the same spherical-vector if and only if they have the same scalar component and the same vector component.
\end{definition}

\begin{definition}\label{def:def2.3}[\textbf{Support of a spherical-vector}]
Let $\alpha$ be a spherical-vector with a non-zero vector component $n$. We define the support of $\alpha$, denoted as $P_{\alpha}$, to be the vector plane for which $n$ serves as a normal vector.
\end{definition}

\begin{remark} \label{rem:rem2.1}
Let $\alpha=(u,v)$ be a spherical-vector with a non-zero vector component $n$.  It is important to note that both vectors $u$ and $v$ belong to the support $P_{\alpha}$ of $\alpha$.
\end{remark}

The following theorem provides a characteristic property of the components of a spherical-vector, which holds significant importance. Indeed, this will allow us to establish an anti-isomorphism between the set of spherical-vectors and the set of unit quaternions in subsequent stages.

\begin{theorem}\label{thm:thm2.1} 
Let $\lambda$ be a real number and $n$ be a vector in $\mathbb{R}^3$. We can assert that $\lambda^2 + \|n\|^2 = 1$ if, and only if, there exists a unique spherical-vector characterized by components $\lambda$ and $n$.
\end{theorem}

\begin{proof}
\begin{enumerate}[$-$]
\item Assume the existence of a spherical-vector $(u, v)$ with components $\lambda$ and $n$. Let $\theta \in [0, \pi]$ represent the non-oriented angle between the two vectors $u$ and $v$. According to \Cref{def:def2.1}, $\lambda = \cos \theta$ and $\|n\| = \sin \theta$. Consequently, $\lambda^2 + \|n\|^2 = 1$.\\

\item Conversely, assume that $\lambda^2 + \|n\|^2 = 1$ and demonstrate the existence of a unique spherical-vector $(u, v)$ with components $\lambda$ and $n$. According to \Cref{def:def2.2}, the uniqueness arises from the condition that two spherical-vectors are equal if and only if they share identical components. It remains now to show the \textit{existence}. To do this, choose any non-zero vector $u$ orthogonal to $n$, and let $v = \lambda u + n \times u$. We begin by verifying that $v$ is non-zero: \[\|v\|^2 = \lambda^2\|u\|^2+\|n\times u\|^2 = \lambda^2\|u\|^2+\|n\|^2 \|u\|^2 = \|u\|^2,\]
where $u$ and $n \times u$ are orthogonal. Consequently, the ordered pair $(u, v)$ represents a spherical-vector with components: \[\frac{u\cdot v}{\|u\| \|v\|}=\frac{u\cdot (\lambda u+n\times u) }{\|u\|^2}=\lambda\] and \[\frac{u \times v}{\|u\| \|v\|} =\frac{u\times (\lambda u+n\times u) }{\|u\|^2}
=\frac{u\times(n\times u)}{\|u\|^2}=n,\]
which can be deduced from the formula, \[u\times(n\times u)=(u\cdot u)n - (u\cdot n)u=\|u\|^2n.\]
\end{enumerate}
\end{proof}

\subsection{Geometric interpretation of a spherical-vector}\label{subsec:subsec2.2}
Let $u$ and $v$ be two non-zero vectors in $\mathbb{R}^3$. According to \Cref{def:def2.1}, the ordered pairs $(u, v)$ and $(u/\|u\|,v/\|v\|)$ share the same components. Consequently, as per \Cref{def:def2.2}, they represent identical spherical-vectors. Given that the vectors in the second ordered pair are unitary, there exists a unique ordered pair $(A, B)$ of points on the unit sphere such that $u/\|u\| = A$, $v/\|v\| = B$, and $(u, v) = (A, B)$. This observation leads to the following remark.

\begin{figure}[H]%
\centering
\includegraphics[width=1\textwidth]{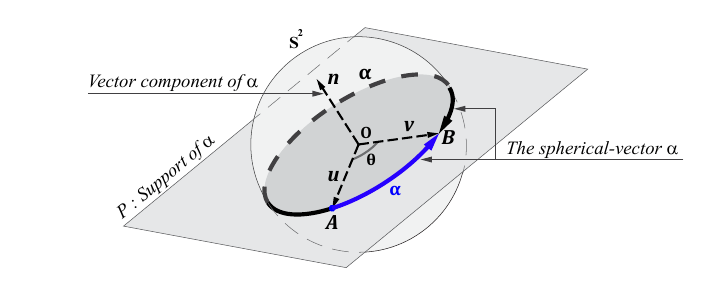}
\caption{Geometric representation of a spherical-vector $\alpha$. The components of $\alpha$ are $\lambda = u \cdot v = \cos \theta$ and $n = u \times v $ where $\|n\| = \sin \theta $ and $\theta$ is the non-oriented angle between the vectors $u$ and $v$.}\label{fig:fig2}
\end{figure}

\begin{remark}\label{rem:rem2.2}
\begin{enumerate}
\item A spherical-vector can be represented by an ordered pair $(u, v)$ of unit vectors in $\mathbb{R}^3$. In this case, based on \Cref{def:def2.1}, the components of $\alpha$ are $\lambda = u \cdot v$ and $n = u \times v$.
\item A spherical-vector can be denoted by an ordered pair $(A, B)$, represented as $\overset{\curvearrowright}{AB}$, consisting of points on the unit sphere $S^2$ (refer to \Cref{fig:fig2}).
\end{enumerate}
\end{remark}

\section{Spherical form of a quaternion}\label{sec:sec3}

In this section, a novel quaternion form is introduced, which facilitates, in \Cref{sec:sec4}, the establishment of a relationship between unit quaternions and spherical-vectors. This relationship unveils intriguing results. For a more comprehensive understanding of quaternion algebra, readers are encouraged to consult (\cite{Kuipers}, Chapter 5) and \cite{Morais}.

\begin{definition} \label{def:def3.1} 
We denote by $\mathbb{V}$ the set of quaternions written in the form $aj + bk + c$, where $(a, b, c) \in \mathbb R^3$.
\end{definition}

\begin{remark} \label{rem:rem3.1} 
Hamilton's approach represents a vector $(a, b, c)$ in $\mathbb{R}^3$ as the pure quaternion $ai + bj + ck$. For further details on this concept, readers are directed to Chapter 5 of \cite{Kuipers} or \cite{Morais}. To prevent confusion with Hamilton's convention, it is important to clarify that the following discussion opts to identify each vector $V = (a, b, c)$ with the quaternion $aj + bk + c \in \mathbb{V}$, which is also denoted by $V$ (with $\mathbb{V}$ defined in \Cref{def:def3.1}). This deliberate choice, which will be elucidated in \Cref{sec:sec5}, allows for the reformulation of the exponential form of quaternions in a manner that adheres to algebraic properties analogous to those of complex numbers' exponential form.
\end{remark}

Let $q = x + yi + zj + tk$, with $(x, y, z, t) \in \mathbb{R}^4$, be a quaternion. We can factorize its imaginary part with respect to $i$ as follows: $q = x + i(tj - zk + y)$. Consequently, we may express $q$ as $x + iV$, where $V$ represents the quaternion $jt - kz + y \in \mathbb{V}$, which is identified with the vector $(t, -z, y)$ according to \Cref{rem:rem3.1}.

\begin{definition}\label{def:def3.2}[Spherical form of a quaternion]
Every quaternion $q$ can be uniquely expressed in the form $x + iV$, where $x \in \mathbb{R}$ and $V \in \mathbb{V}$. This representation is referred to as the \textit{spherical form} of $q$. The \textit{spherical components} of $q$ consist of the real number $x$ (called the \textit{scalar component}) and the vector $V$ (referred to as the \textit{vector component}).
\end{definition}

\begin{proposition} \label{propo:propo3.1} 
\begin{enumerate}
\item Two quaternions are equal if and only if their spherical components are the same.
\item Given a quaternion $q$ with spherical form $x + iV$, the squared magnitude of $q$ is given by $\vert q\vert^2 = x^2 + \|V\|^2$.
\end{enumerate}
\end{proposition}
\begin{proof}
immediate.
\end{proof}

The subsequent lemma presents several valuable properties of the spherical form, which will be utilized later in this article.

\begin{lemma}\label{lem:lem3.1}
For all vectors $u$ and $v$, the following relationships holds.
\begin{enumerate}
\item $iu = \overline{u} i$,
\item $uv = v \cdot \overline{u} + i(v \times \overline{u})$,
\item If $v$ is a unit vector then $v^{-1}u = u \cdot v + i (u\times v) $.
\end{enumerate}
\end{lemma}

\begin{proof}
\begin{enumerate}
\item Let $u = (x, y, z)$ be a vector in $\mathbb{R}^3$. In accordance with \Cref{rem:rem3.1}, we identify $u$ with the quaternion $jx + ky + z$. Thus, we obtain \[iu = i(jx + ky + z) = kx - jy + iz = (-jx - ky + z)i = \overline{u}i.\]

\item Let $u=(x,y,z)$ and $v=(x',y',z')$ be two vectors in $\mathbb R^3$. We have
\begin{align*}
uv &=(jx+ky+z)(jx'+ky'+z') \\
&=(-xx'-yy'+zz')+i\left[j(yz'+y'z)-k(xz'+x'z)+(xy'-x'y)\right] \\
&=(-xx'-yy'+zz')+i \begin{pmatrix} yz'+y'z \\ -xz'-x'z \\ xy'-x'y \end{pmatrix} \\
&=v\cdot \overline{u}+i(v\times \overline{u}).
\end{align*}
\item The result follows from (ii) and the fact that the inverse of unit vector is equal to its conjugate.
\end{enumerate}
\end{proof}

\section{Algebraic structure of spherical-vectors}\label{sec:sec4}
Throughout the subsequent discussion, let $\mathcal{S}$ denote the set of spherical-vectors and $\mathbb{H}_1$ represent the multiplicative group of unit quaternions. In this section, we define and examine an additive group structure on $\mathcal{S}$, facilitated by a bijective map from $\mathcal{S}$ to $\mathbb{H}_1$.

\subsection{Additive group of spherical-vectors}\label{subsec:subsec4.1}

Let $\alpha$ be a spherical-vector with components $\lambda$ and $n$, denoted as $\alpha(\lambda,n)$. It is straightforward to associate $\alpha$ with the quaternion $q_{\alpha}$, defined by its spherical form $\lambda + in$. According to \Cref{thm:thm2.1} and \Cref{propo:propo3.1}.(ii), the quaternion $q_{\alpha}$ is unitary.

\begin{theorem}\label{thm:thm4.1} 
The mapping $\mu $ that associates each spherical-vector $\alpha(\lambda,n)$ with the unit quaternion $\lambda + in$ is bijective.
\begin{equation}\label{eqt:eqt4.1}
\mu :\ \mathcal{S}\ \longrightarrow\ \mathbb{H}_1 , \ \alpha(\lambda,n)\ \longmapsto\ \lambda+in. 
\end{equation}
\end{theorem}

\begin{proof}
Let $q = x + iw$ be a unit quaternion expressed in its spherical form. According to \Cref{thm:thm2.1}, there exists a unique spherical-vector $\alpha$ with components $x$ and $w$. Thus, $\mu (\alpha) = q$.
\end{proof}

\begin{proposition}\label{propo:propo4.1} 
Let $\alpha = (u, v)$ be a spherical-vector represented by an ordered pair of unit vectors. We obtain 
\begin{equation} \label{eqt:eqt4.2}
\mu({\alpha}) = u \cdot v + i (u \times v) = v^{- 1} u.
\end{equation}
\end{proposition}

\begin{proof}
The first expression of $\mu(\alpha)$ is derived from the definition of $\mu$ in (\ref{eqt:eqt4.1}) and the fact that $u \cdot v$ and $u \times v$ are the components of $\alpha$ according to \Cref{def:def2.1}. The second expression is a consequence of \Cref{lem:lem3.1}.(iii).
\end{proof}

The mapping $\mu$ enables us to transfer the structure of $\mathbb{H}_1$ to $\mathcal{S}$. We can then define an additive internal composition law on $\mathcal{S}$ as $\alpha + \beta = \mu^{-1}(\mu(\beta) \mu(\alpha))$, which implies that $\mu(\alpha + \beta) = \mu(\beta) \mu(\alpha)$. Consequently, $\mu$ serves as an \textit{anti-isomorphism} of groups, leading to the following theorem.

\begin{theorem}\label{thm:thm4.2} 
The mapping $\mu$ is an \textit{anti-isomorphism} of groups, and the set $\mathcal{S}$ of spherical-vectors is a non-commutative additive group isomorphic to the group $\mathbb{H}_1$ of unit quaternions.
\end{theorem}

\subsection{Algebraic properties of spherical-vectors}\label{subsec:subsec4.2}

In this section, using the anti-isomorphism $\mu$ established in \Cref{subsec:subsec4.1}, we transfer the algebraic properties of the multiplicative group $\mathbb {H}_1$  of unit quaternions into the additive group $\mathcal {S}$ of spherical-vectors.

\subsubsection{Zero spherical-vector}\label{subsubsec:subsubsec4.2.1}

\begin{proposition}\label{propo:propo4.2} 
The neutral element of $\mathcal{S}$ which we call zero spherical-vector and denote by $\overset{\curvearrowright}{0}$, is represented by any ordered pair $(u, u)$ where $u$ is a unit vector.
\end{proposition}
\begin{proof}
Suppose the zero spherical-vector $\overset{\curvearrowright}{0}$ is represented by an ordered pair $(u,v)$ of unit vectors. Given that $\mu$ is a group anti-isomorphism, we have $\mu\left((u,v)\right) = 1$. According to (\ref{eqt:eqt4.2}), $\mu\left((u,v)\right) = v^{-1}u$. Consequently, $\mu\left((u,v)\right) = 1$ holds true if and only if $u=v$.
\end{proof}

\subsubsection{Opposite of a spherical-vector}\label{subsubsec:subsubsec4.2.2}

\begin{proposition}\label{propo:propo4.3} 
\begin{enumerate}
\item Let $u$ and $v$ be two unit vectors. The opposite of the spherical-vector $(u,v)$, which we denote by $-(u, v)$, is the spherical-vector $(v,u)$.
\item If $\lambda$ and $n$ are the components of a spherical-vector $\alpha$, then the  components of $-\alpha$ are $\lambda$ and $-n$.
\end{enumerate} 
\end{proposition}

\begin{proof}
\begin{enumerate}
\item Given that $\mu$ represents an anti-isomorphism of groups, according to (\ref{eqt:eqt4.2}), the proof stems from the identity \[\mu\left(-(u,v)\right) = \left(\mu \left ((u, v) \right) \right)^{-1} = (v^{-1}u)^{-1} = u^{-1}v = \mu\left((v,u)\right).\]
\item This is an easy consequence of \Cref{def:def2.1}.
\end{enumerate} 
\end{proof}

\subsubsection{Chasles relation}\label{subsubsec:subsubsec4.2.4}

\begin{proposition}\label{propo:propo4.4} 
Let $u$, $v$ and $w$ be three unit vectors. We have \[(u,v)+(v,w)=(u,w).\]
\end{proposition}
\begin{proof}
As $\mu$ represents an anti-isomorphism of groups and according to (\ref{eqt:eqt4.2}), the Chasles relation can be readily derived from the following equation \[\mu\left((u,v)+(v,w)\right)=\mu\left((v,w)\right)\mu\left((u,v)\right)=w^{-1}u=\mu\left((u,w)\right).\]
\end{proof}

\subsubsection{Subtraction operation on $\mathcal{S}$}\label{subsubsec:subsubsec4.2.3}
\begin{definition}\label{def:def4.1} 
Let $\alpha$ and $\beta$ be two spherical-vectors.
We define the subtraction operation on $\mathcal{S}$ by, $\alpha - \beta = \alpha + (-\beta)$.
\end{definition}

\subsubsection{Straight spherical-vector}\label{subsubsec:subsubsec4.2.5}

Let $u$ be a unit vector. According to \Cref{def:def2.1}, the components of the spherical-vector $(u, -u)$ are $-1$ and $\vec{0}$, and they are independent of $u$. This observation leads to the following definition.

\begin{definition}\label{def:def4.2} 
We define the \textit{straight spherical-vector}, denoted by $\overset{\curvearrowright}{\pi}$, as the spherical-vector with components $-1$ and $\vec{0}$, represented by any ordered pair $(u, -u)$, where $u$ is a unit vector.
\end{definition}

Let $\alpha$ be a spherical-vector represented by an ordered pair $(u, v)$ of unit vectors. According to \Cref{def:def2.1}, the vector component of $\alpha$ is $n_{\alpha} = u \times v$. Given that $u$ and $v$ are unit vectors, the vector $n_{\alpha}$ is zero if and only if $v = u$ or $v = -u$. Consequently, by invoking \Cref{propo:propo4.2} and \Cref{def:def4.2}, we arrive at the following remark.

\begin{remark}\label{rem:rem4.1} 
The vector component of a spherical-vector is zero if and only if this spherical-vector is zero or straight.
\end{remark}

\subsection{Sum of two spherical-vectors}\label{subsec:subsec4.3}

In this section we return with more details to the sum of two spherical-vectors defined by their components. To do so, we will need the following lemmas.

\begin{figure}[H]%
\centering
\includegraphics[width=1\textwidth]{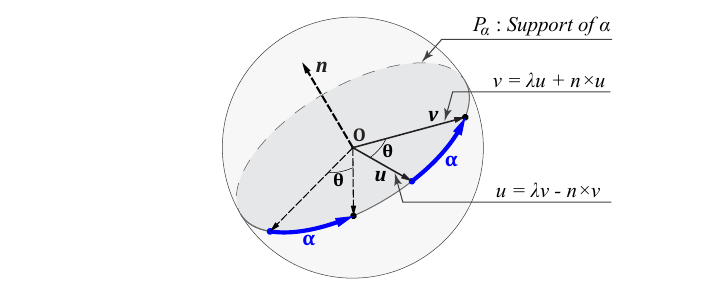}
\caption{Spherical-vector $\alpha(\lambda, n)$ can be decomposed into a pair of unit vectors $(u,v)$, revealing the relationship between the spherical components of $\alpha$ and its corresponding unit vectors.}
\label{fig:fig3}
\end{figure}

\begin{lemma}\label{lem:lem4.1} 
Let $\alpha$ be a spherical-vector with components $\lambda$ and $n$, where $n \neq \vec{0}$. Let $u$, $v$ and $w$ be three unit vectors $(\Cref{fig:fig3})$. We have
\begin{enumerate}
\item The equality $(u,v) = \alpha$ is satisfied if and only if $v = u(\lambda-in)$, which holds if and only if $u$ is orthogonal to $n$ and $v$ can be expressed as $v = \lambda u + n \times u$.
\item The equality $(u,v) = \alpha$ is satisfied if and only if $u = v (\lambda + in)$, which holds if and only if $v$ is orthogonal to $n$ and $u = \lambda v - n \times v$.
\end{enumerate}
\end{lemma}

\begin{proof}
\begin{enumerate}
\item Based on (\ref{eqt:eqt4.1}) and (\ref{eqt:eqt4.2}), we have the following equivalence:
\[(u,v) = \alpha \Leftrightarrow \mu((u,v)) = \mu(\alpha) \Leftrightarrow v^{-1}u = \lambda + in\Leftrightarrow v = u(\lambda - in).\]
According to \Cref{rem:rem2.1}, if $(u,v) = \alpha$, it follows that $u \perp n$. Therefore,
\[(u,v) = \alpha \Leftrightarrow u \perp n \text{ and } v = u(\lambda - in).\]
Finally, using \Cref{lem:lem3.1}.(i) and \Cref{lem:lem3.1}.(ii), when $u \perp n$, we can deduce that \[u(\lambda - in) = \lambda u - i\overline{u}n = \lambda u + n \times u.\] Thus, the result is obtained.
\item This outcome can be readily derived from part (i) by noting that, according to \Cref{propo:propo4.3}, $(v,u) = -(u,v)$ and the components of $-\alpha$ are $\lambda$ and $-n$.
\end{enumerate}
\end{proof}

\begin{figure}[H]%
\centering
\includegraphics[width=1\textwidth]{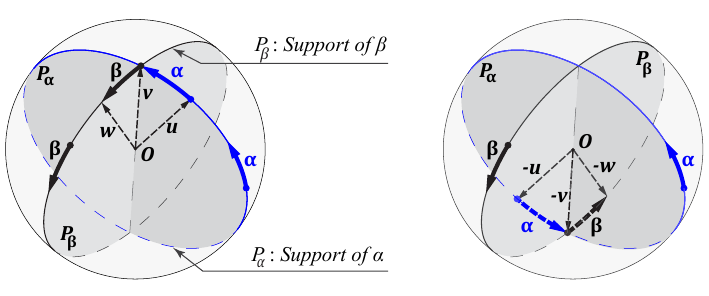}
\caption{$n_{\alpha} \times n_{\beta} \neq \vec{0}$ : The supports of $\alpha$ and $\beta$ intersect in a vector line spanned by the vector $n_{\alpha} \times n_{\beta}$.}\label{fig:fig4}
\end{figure}

\begin{figure}[H]%
\centering
\includegraphics[width=1\textwidth]{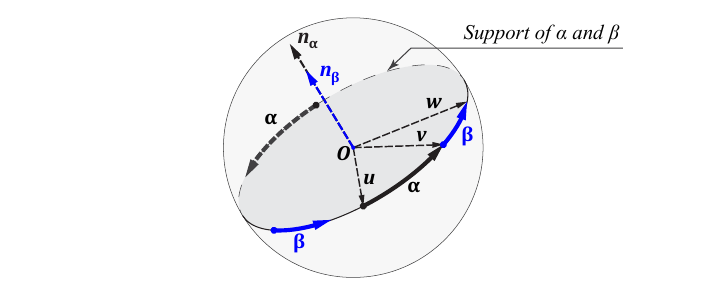}
\caption{$n_{\alpha} \times n_{\beta} = \vec{0}$ : The two spherical-vectors $\alpha$ and $\beta$ have the same support.}\label{fig:fig5}
\end{figure}

\begin{lemma} \label{lem:lem4.2} 
Let $\alpha$ and $\beta$ be two spherical-vectors with respective supports $P_{\alpha}$ and $P_{\beta}$, and components $(\lambda_{\alpha}, n_{\alpha})$ and $(\lambda_{\beta}, n_{\beta})$, where $n_{\alpha} \neq \vec{0}$ and $n_{\beta} \neq \vec{0}$. We proceed by considering two distinct cases.
\begin{enumerate}
\item In the case where $n_{\alpha} \times n_{\beta} \neq \vec{0}$, there exist precisely two sets of three unit vectors $u$, $v$, and $w$ satisfying the conditions $\alpha=(u,v)$ and $\beta=(v,w)$, as illustrated in \Cref{fig:fig4}. These sets are as follows,
\[\begin{cases} v=\frac{n_{\alpha} \times n_{\beta}}{\Vert n_{\alpha} \times n_{\beta}\Vert} \\ u=\lambda_{\alpha}v-n_{\alpha}\times v\\ w=\lambda_{\beta}v+n_{\beta}\times v \end{cases} \text{and}\quad
\begin{cases} v'=-v \\ u'=-u\\ w'=-w. \end{cases}\]
\item If $n_{\alpha} \times n_{\beta} = \vec {0}$, the two spherical-vectors $\alpha$ and $\beta$ have the same support $P$. In this case, for any unit vector $v$ in $P$, there exist exactly two unique vectors $u=\lambda_{\alpha} v - n_{\alpha} \times v$ and $w=\lambda_{\beta} v + n_{\beta} \times v$ that satisfy the conditions $\alpha=(u,v)$ and $\beta=(v,w)$, as illustrated in \Cref{fig:fig5}.
\end{enumerate}
\end{lemma}

\begin{figure}[H]%
\centering
\includegraphics[width=1\textwidth]{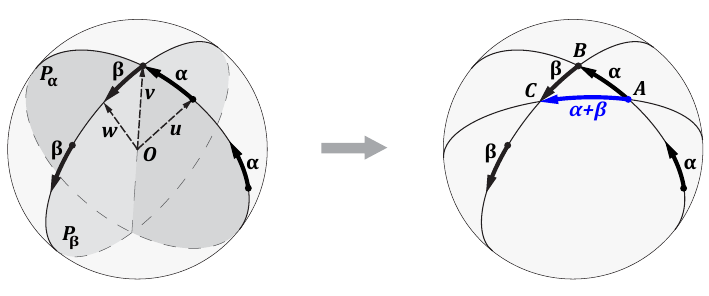}
\caption{Addition of spherical-vectors on $S^{2}$: $\overset{\curvearrowright}{AB}+\overset{\curvearrowright}{BC}=\overset{\curvearrowright}{AC}.$}
\label{fig:fig6}
\end{figure}

\begin{proof}
\begin{enumerate}
\item The objective is to determine the solution to the following system of equations, the unknowns being the unit vectors $u$, $v$ and $w$.\[\begin{cases}
(u,v)=\alpha\\
(v,w)=\beta,
\end{cases}\]
As stated in \Cref{lem:lem4.1}, we can rewrite this system as follows:
\[\begin{cases} v \perp n_{\alpha} \text{ and } u = \lambda_{\alpha} v - n_{\alpha} \times v \\ v \perp n_{\beta} \text{ and } w = \lambda_{\beta} v + n_{\beta} \times v.\end{cases}\]
If we assume that $v$ is orthogonal to both $n_{\alpha}$ and $n_{\beta}$, and since $v$ is a unit vector, we can write $v$ as $v=\pm\frac{n_{\alpha} \times n_{\beta}}{\Vert n_{\alpha} \times n_{\beta}\Vert}$. Consequently, we can proceed to solve the system with the knowledge of $v$, and then use the expressions for $u$ and $w$ given above. As a result, we are left with only two possible solutions: \[\begin{cases} v=\frac{n_{\alpha} \times n_{\beta}}{\Vert n_{\alpha} \times n_{\beta}\Vert} \\ u = \lambda_{\alpha} v - n_{\alpha} \times v \\ w = \lambda_{\beta} v + n_{\beta} \times v\\ \end{cases} \text{or} \quad
\begin{cases} v=-\frac{n_{\alpha} \times n_{\beta}}{\Vert n_{\alpha} \times n_{\beta}\Vert} \\ u = \lambda_{\alpha} v - n_{\alpha} \times v \\ w = \lambda_{\beta} v + n_{\beta} \times v. \end{cases}\]\\

\item This is an immediate consequence of \Cref{lem:lem4.1}.
\end{enumerate}
\end{proof}
We now possess the essential tools to construct the sum of two spherical-vectors $\alpha$ and $\beta$, with respective vector components $n_{\alpha}$ and $n_{\beta}$. Two cases must be considered:

\begin{enumerate}[$-$]
\item If both $n_{\alpha}$ and $n_{\beta}$ are non-zero (refer to \Cref{fig:fig6}), then according to \Cref{lem:lem4.2}, there exist three unit vectors $u$, $v$, and $w$ such that $\alpha = (u, v)$ and $\beta = (v, w)$. Therefore, based on the Chasles relation in \Cref{propo:propo4.4}, we obtain: \[\alpha + \beta = (u, v) + (v, w) = (u, w).\]
\item If either $n_{\alpha}$ or $n_{\beta}$ is zero, let's assume $n_{\beta}$ (the other case follows a similar logic), then by \Cref{rem:rem4.1}, either $\beta = \overset{\curvearrowright}{0}$ or $\beta = \overset{\curvearrowright}{\pi}$. Let $u$ and $v$ be two unit vectors such that $\alpha = (u,v)$. If $\beta = \overset{\curvearrowright}{0}$, then $\alpha + \beta = \alpha$. If $\beta = \overset{\curvearrowright}{\pi}$, then according to \Cref{def:def4.2}, the straight spherical-vector $\overset{\curvearrowright}{\pi}$ can be represented by $(v,-v)$. Consequently, \[\alpha + \beta = (u,v) + (v,-v) = (u,-v).\]
\end{enumerate}

\section{Argument and polar form of a quaternion}\label{sec:sec5}
\begin{figure}[H]%
\centering
\includegraphics[width=1\textwidth]{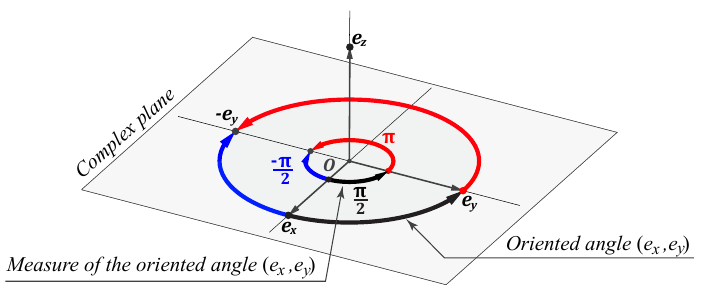}
\caption{Spherical-vectors generalize oriented angles between two vectors to $\mathbb{R}^3$ space.}
\label{fig:fig7}
\end{figure}

In this section, we reformulate the definition of the argument of a quaternion through the concept of spherical-vectors, enabling us to rewrite the polar form of a quaternion and establish all algebraic properties of the exponential representation (refer to Chapter 5 of \cite{Kuipers}, \cite{Lam}, and \cite{Quaternions} for further details on the polar and exponential forms of quaternions). In the oriented complex plane, each oriented angle $(u,v)$ between two vectors $u$ and $v$ is identified with its measure $\theta$ in $]-\pi,\pi]$. Consequently, the fundamental algebraic property of the exponential form of unit complex numbers can be interpreted in two equivalent ways: either by utilizing the geometric aspect of oriented angles between two vectors (i.e., the spherical-vectors they represent) or by their numerical aspect (i.e., their measures). For instance, we have $i=e^{i(e_x,e_y)} = e^{i\pi/2}$ and $-1=e^{i(e_y,-e_y)} = e^{i \pi}$ (refer to \Cref{fig:fig7}). Therefore, the product $i\times(-1)$ can be represented by the following expressions (See \cite{Andreescu}, Chapter 2) \[i\times(-1) = e^{i(e_x,e_y)} \times e^{i(e_y,-e_y)} = e^{i\left((e_x,e_y)+(e_y,-e_y)\right)} = e^{i(e_x,-e_y)}\] or \[i\times(-1) = e^{i\frac{\pi}{2}} \times e^{i\pi} = e^{i\left(\frac{\pi}{2} + \pi\right)} = e^{i\frac{3\pi}{2}} = e^{-i\frac{\pi}{2}}.\]

Later on, we will demonstrate that the first formula presented above, involving $i \times (-1)$, can be naturally extended to the set $\mathbb{H}_1$ of unit quaternions.

\subsection{Argument of a quaternion}\label{subsec:subsec5.1}

\begin{definition}\label{def:def5.1}
Let $q$ be a non-zero quaternion. According to \Cref{thm:thm4.1}, the mapping $\mu$ is bijective from $\mathcal{S}$ to $\mathbb{H}_1$. We define the \textit{argument} of $q$ as the spherical-vector given by:
\begin{equation}\label{eqt:eqt5.1}
\arg(q)=\mu^{-1}\left(\frac{q}{ \vert q \vert}\right).
\end{equation}
In other words, given a spherical-vector $\alpha$,
\begin{equation}\label{eqt:eqt5.2}
\arg(q)=\alpha \text{ if and only if } \mu(\alpha)=\frac{q}{\vert q \vert}.
\end{equation}
\end{definition}

\begin{remark}\label{rem:rem5.1}
Let $q$ be a unit quaternion and $\alpha$ a spherical-vector whose components are $\lambda$ and $n$.
\begin{enumerate}
\item By (\ref{eqt:eqt5.2}) and (\ref{eqt:eqt4.1}), $\arg(q)=\alpha \text{ if and only if } q = \mu(\alpha) = \lambda+in$. Therefore, the spherical components of $q$ and the components of its argument $\alpha$ are the same.
\item If $\alpha$ is represented by an ordered pair $(u,v)$ of unit vectors, then by (i) and (\ref{eqt:eqt4.2}), $\arg(q)=\alpha \text{ if and only if } q = u.v + i (u \times v)$.
\end{enumerate}
\end{remark} 

\subsection{Cosine and sine of a spherical-vector}\label{subsec:subsec5.2}

We have seen in this paper that spherical-vectors are order pairs of non-zero vectors in the oriented space $\mathbb{R}^3$. Consequently, the oriented angles between two vectors in the oriented plane $\mathbb{R}^2$, represent particular instances of spherical-vectors. We would therefore like to define mappings that extend the cosine and sine functions of oriented angles to the set $\mathcal{S}$ of spherical-vectors. We illustrate the idea behind this extension through the subsequent example.

\begin{figure}[H]%
\centering
\includegraphics[width=1\textwidth]{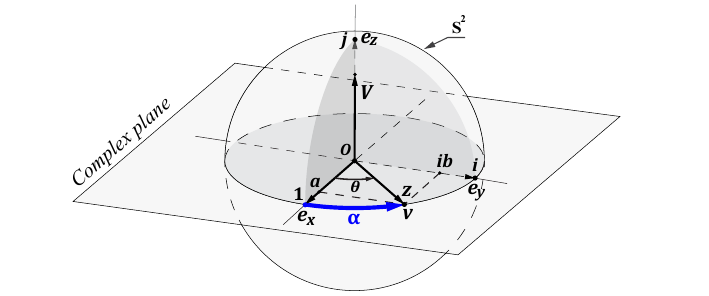}
\caption{The spherical-vector $\alpha=(e_x,v)$ is an oriented angle between two vectors in the complex plane, and $\theta$ is the measure of the oriented angle $(e_x,v)$.}
\label{fig:fig8}
\end{figure}

Let $z = a + ib$ denote a unit complex number, represented by a unit vector $v(a,b)$ in the complex plane, with $a,b \in \mathbb{R}$. The argument of $z$ can be expressed either as the oriented angle $\alpha = (e_x, v)$, which constitutes a spherical-vector, or as its measure $\theta \in ]-\pi, \pi]$ (refer to Figure \ref{fig:fig8}). Consequently, we arrive at the following equation:
\begin{equation}\label{eqt:eqt5.3}
z = \cos\theta + i \sin\theta = \cos\alpha + i \sin\alpha.
\end{equation}
On the other hand, according to \Cref{def:def3.2}, the expression $a + ib$ represents the spherical form of $z$. This can be written as:
\begin{equation}\label{eqt:eqt5.4}
z = a + i(j \times 0 + k \times 0 + b) = a + i(0, 0, b) = a + iV,
\end{equation}

where $V=(0, 0, b)=b e_z$ denotes the vector component of $z$ associated with the real number $b \in \mathbb{V}$, as per \Cref{rem:rem3.1}. By employing Equations (\ref{eqt:eqt5.3}) and (\ref{eqt:eqt5.4}), we obtain:
\[\begin{cases}
\cos\alpha = a \in \mathbb{R} \\
\sin\alpha = V \in \mathbb{V},
\end{cases}\]
where $a$ and $V$ also serve as the components of $\alpha$, according to \Cref{rem:rem5.1}(i).
This leads us to the next definition.

\begin{definition} \label{def:def5.2}
Let $\alpha$ represent a spherical-vector with components $\lambda$ and $n$.
\begin{enumerate}[$-$]
\item The \textit{cosine} of $\alpha$, denoted as $\cos\alpha$, refers to the scalar component $\lambda$.
\item The \textit{sine} of $\alpha$, denoted as $\sin\alpha$, refers to the vector component $n$, which is identified with a quaternion in $\mathbb{V}$ according to \Cref{rem:rem3.1}.
\end{enumerate}
Consequently, we establish two mappings on $\mathcal{S}$ as follows:
\[\begin{array}{r@{\ }l}
\cos :\ \mathcal{S} & \longrightarrow\ \mathbb{R} \\
\alpha & \longmapsto\ \lambda \\
\end{array} \qquad;\quad  \begin{array}{r@{\ }l}
\sin:\ \mathcal{S}  & \longrightarrow\ \mathbb{V} \\
\alpha & \longmapsto\ n. \\
\end{array}\]
\end{definition}

\begin{remark}\label{rem:rem5.2}
Let $\alpha$ be a spherical-vector represented by an ordered pair $(u, v)$ of unit vectors.
\begin{enumerate}
\item By \Cref{def:def2.1}, the components of $\alpha$ are $u \cdot v$ and $u\times v$. So, by \Cref{def:def5.2} above,
\begin{equation}\label{eqt:eqt5.5}
\begin{cases}
\cos\alpha=u \cdot v \\
\sin\alpha=u\times v \in \mathbb{V}.
\end{cases}
\end{equation}
Thus by (\ref{eqt:eqt4.2}),
\begin{equation}\label{eqt:eqt5.6}
\mu(\alpha) = \cos\alpha + i \sin\alpha.
\end{equation}

\item Let $\theta$ denote the non-oriented angle between the two vectors $u$ and $v$, where $\theta \in [0, \pi]$. Based on (\ref{eqt:eqt5.5}), we observe that $\cos\alpha = \cos\theta$ and $|\sin\alpha| = |u \times v| = \sin\theta$ (\Cref{fig:fig2}). Consequently, the following relationship is established:
\begin{equation}\label{eqt:eqt5.7}
\cos^2\alpha+\vert \sin\alpha \vert^2=1,
\end{equation}
\end{enumerate}
\end{remark}

\begin{example}\label{exm:exm5.1}
Consider the unit vectors $u = \left(\frac {\sqrt {2}} {2}, \frac {\sqrt{2}}{2}, 0\right)$ and $v = \left(\frac {\sqrt {3}}{3}, \frac {\sqrt{3}} {3}, \frac {\sqrt{3}}{3}\right)$. The components of the spherical-vector $\alpha$, represented by the ordered pair $(u, v)$, can be expressed as:
\[\begin{cases}
\lambda = u \cdot v = \frac {\sqrt {6}}{3}\\
n = u \times v = \left(\frac {\sqrt {6}}{6}, - \frac {\sqrt {6}}{6}, 0\right).
\end{cases}\]
According to \Cref{rem:rem3.1}, $n=\frac{\sqrt{6}}{6}j-\frac{\sqrt{6}}{6}k\in \mathbb{V}$. Therefore, based on (\ref{eqt:eqt5.5})
\[\begin{cases}
\cos\alpha=\frac{\sqrt{6}}{3}\\
\sin\alpha=\frac{\sqrt{6}}{6}j-\frac{\sqrt{6}}{6}k.
\end{cases}\]
\end{example} 

\subsection{Polar form of a quaternion}\label{subsec:subsec5.3}

\begin{theorem}\label{thm:thm5.1}
Let $q$ denote a non-zero quaternion, $r>0$ represent a positive real number, and $\alpha$ be a spherical-vector. The equality $q = r (\cos\alpha + i \sin\alpha)$ holds if and only if $\vert q \vert = r$ and $\arg(q) = \alpha$. The expression $r (\cos\alpha + i \sin\alpha)$ is referred to as the \textit{polar form} of $q$.
\end{theorem}
\begin{proof}
\begin{enumerate}[$-$]
\item Assume that $q = r(\cos\alpha + i\sin\alpha)$. According to \Cref{propo:propo3.1}.(ii) and (\ref{eqt:eqt5.7}), we derive that $\vert q\vert^2 = r^2(\cos^2(\alpha) + \vert \sin(\alpha) \vert^2) = r^2$. Furthermore, based on (\ref{eqt:eqt5.6}), we have $\mu(\alpha) = \frac{q}{\vert q \vert}$. Consequently, by referring to  (\ref{eqt:eqt5.2}), it follows that $\arg(q) = \alpha$.\\
\item Conversely, assume that $\vert q \vert = r$ and $\arg(q) = \alpha$. By (\ref{eqt:eqt5.2}) and (\ref{eqt:eqt5.6}), we obtain the relation $q/r = \mu(\alpha) = \cos(\alpha) + i\sin(\alpha)$.
\end{enumerate}
\end{proof}

\subsection{Exponential form of a quaternion}\label{subsec:subsec5.4}

\begin{definition} \label{def:def5.3}
Let $\alpha$ denote a spherical-vector. The expression $e^{i \alpha}$ refers to the unit quaternion corresponding to the argument $\alpha$. Based on \Cref{thm:thm5.1}, we can express this relationship as:
\begin{equation}\label{eqt:eqt5.8}
e^{i \alpha} = \cos\alpha + i \sin \alpha.
\end{equation}
Hence, if $q$ is a quaternion with a modulus $r > 0$ and an argument $\alpha$, by applying \Cref{thm:thm5.1} and (\ref{eqt:eqt5.8}), we derive:
\begin{equation}\label{eqt:eqt5.9}
q = r (\cos\alpha + i \sin\alpha) = r e^{i\alpha}.
\end{equation}
\end{definition}

\begin{remark}
Assuming $\alpha$ is a spherical-vector, equations (\ref{eqt:eqt5.6}) and (\ref{eqt:eqt5.8}) yield the expression
\begin{equation}\label{eqt:eqt5.10}
\mu(\alpha) = e^{i\alpha}.
\end{equation}
\end{remark} 

\subsection{Algebraic properties of the exponential form}\label{subsec:subsec5.5}

\begin{proposition} \label{propo:propo5.1}
Let $\alpha$ and $\beta$ represent two spherical-vectors, and let $m$ be an integer. The following properties hold:
\begin{align*}
&(1)\ e^{i(\alpha+\beta)}=e^{i\beta}e^{i\alpha} &&(2)\ \left(e^{i\alpha}\right)^{-1}=\overline{e^{i\alpha}}=e^{-i\alpha} \\
&(3)\ e^{i(\alpha-\beta)}=\left(e^{i\beta}\right)^{-1}e^{i\alpha}&&(4)\ e^{i\overset{\curvearrowright}{\pi}}=-1 \\
&(5)\ -e^{i\alpha}=e^{i(\overset{\curvearrowright}{\pi}+\alpha)}=e^{i(\alpha+\overset{\curvearrowright}{\pi})} &&(6)\ \left(e^{i\alpha}\right)^{m}=e^{im\alpha }. 
\end{align*}
\end{proposition}

\begin{proof}
\begin{enumerate}[$(1)$]
\item Based on \Cref{thm:thm4.2}, $\mu$ is an anti-isomorphism of groups. Consequently, $\mu(\alpha + \beta) = \mu(\beta) \mu(\alpha)$. According to  (\ref{eqt:eqt5.10}), we derive that $e^{i (\alpha + \beta)} = e^{i\beta} e^{i\alpha}$.

\item Similarly, we have $\mu(-\alpha) = \mu(\alpha)^{-1}$. Therefore, by (\ref{eqt:eqt5.10}), $e^{-i\alpha} = \left(e^{i\alpha}\right)^{-1}$. As $e^{i\alpha}$ is a unit quaternion, it follows that $\overline{e^{i\alpha}} = \left(e^{i\alpha} \right)^{-1}$.

\item To establish this property, one can simply combine $(1)$ and $(2)$.

\item According to \Cref{def:def4.2}, the components of $\overset{\curvearrowright}{\pi}$ are $-1$ and $\vec{0}$. Thus, by (\ref{eqt:eqt4.1}) and (\ref{eqt:eqt5.10}), we obtain $\mu(\overset{\curvearrowright}{\pi}) = -1 = e^{i\overset{\curvearrowright}{\pi}}$.

\item To prove this property, one can combine $(1)$ and $(4)$.

\item This result can be easily demonstrated through recursion on $m \in \mathbb{N}$, using $(1)$. To extend it to $\mathbb{Z}$, we apply $(2)$.
\end{enumerate}
\end{proof}

We can restate the aforementioned properties as follows. 

\begin{proposition}\label{propo:propo5.2}
Let $p$ and $q$ be two non-zero quaternions, and let $m$ be an integer. The following properties hold.
\begin{align*}
&\bullet\ \arg(pq)=\arg(q)+\arg(p) &&\bullet\ \arg(q^{-1})=\arg(\overline{q})=-\arg(q) \\
&\bullet\ \arg(p^{-1}q)=\arg(q)-\arg(p) && \bullet\ \arg(-1)=\overset{\curvearrowright}{\pi} \\
&\bullet\ \arg(-q)=\overset{\curvearrowright}{\pi}+\arg(q)=\arg(q)+\overset{\curvearrowright}{\pi} &&
\bullet\ \arg(q^{m})=m\arg(q).
\end{align*}
\end{proposition}

\section{Examples of applications}\label{sec:sec6}
As observed in previous sections, the anti-isomorphism $\mu$ enables us to associate every unit quaternion $q$ with its argument, which is a spherical-vector. Consequently, we can conveniently represent unit quaternions and their multiplication geometrically through spherical-vectors.

\subsection{Practical method for determining the argument of a unit quaternion}\label{subsec:subsec6.1}
Let $q = \lambda + in$ be a unit quaternion expressed in its spherical form, with $\alpha$ denoting its argument. According to \Cref{rem:rem5.1}.(i), the components of $\alpha$ are $\lambda$ and $n$. The approach discussed here involves representing the spherical-vector $\alpha$ as an ordered pair $(u,v)$ of unit vectors.

\begin{enumerate}[$-$]
\item Assume that $n = \vec{0}$. According to \Cref{rem:rem4.1}, $\alpha$ can be represented as $\overset{\curvearrowright}{0} = (u,u)$ or $\alpha = \overset{\curvearrowright}{\pi} = (u, -u)$, where $u$ is an arbitrary unit vector.
\item Now, let us consider the case where $n \neq \vec{0}$. According to \Cref{lem:lem4.1}.(i),
\begin{equation} \label{eqt:eqt6.1}
(u,v) = \alpha\text{ if and only if }\left(u \perp n\text{ and }v = \lambda u + n \times u = u\overline{q}\right).
\end{equation}
To determine an ordered pair $(u,v)$ of unit vectors representing the argument $\alpha$ of $q$, we simply need to select any unit vector $u$ orthogonal to $n$, and then determine $v$ using one of the two formulas $v = u\overline{q} = \lambda u + n \times u$.
\end{enumerate}

\begin{example}[\textbf{Argument of quaternions $\boldsymbol{i}$, $\boldsymbol{j}$ and $\boldsymbol{k}$}]\text{} \label{exm:exm6.1}
\begin{enumerate}[$\bullet$]
\item \textbf{Argument of $\boldsymbol{i.}$} Using \Cref{def:def3.2} and \Cref{rem:rem3.1}, the spherical form of $i$ is given by \[i = i (0j + 0k + 1) = i (0,0,1) = ie_z.\]
Consequently, the spherical components of $i$ are $x=0$ and $w=e_z$. From (\ref{eqt:eqt6.1}), to represent the argument of $i$ as an ordered pair $(u,v)$ of unit vectors, we observe that $e_x$ is orthogonal to $e_z$. Therefore, we can select $u=e_x$ and $v=e_z \times e_x = e_y$. As a result, the argument of $i$ corresponds to the spherical-vector $\alpha_i = (e_x, e_y)$ (\Cref{fig:fig9}). In accordance with \Cref{def:def5.3},
\begin{equation} \label{eqt:eqt6.2}
i = e^{i(e_x, e_y)} \ \text{with} \ \arg(i) = (e_x, e_y).
\end{equation}

\item \textbf{Argument of $\boldsymbol{j.}$} Given that \[j=-ik=-i(0j+1k+0)=-i(0,1,0)=-ie_y,\] the spherical components of $j$ are $0$ and $-e_y$. An ordered pair $(u, v)$ of unit vectors represents the argument of $j$ if and only if $u \perp e_y$ and $v = -e_y \times u$. As $e_x \perp e_y$, we can conveniently choose $u = e_x$ and $v = -e_y \times e_x = e_z$. Consequently, the argument of $j$ corresponds to the spherical-vector $\alpha_j = (e_x, e_z)$ (\Cref{fig:fig9}). Thus,
\begin{equation} \label{eqt:eqt6.3}
j = e^{i(e_x, e_z)} \ \text{with} \ \arg(j) = (e_x, e_z).
\end{equation} 

\item \textbf{Argument of $\boldsymbol{k.}$} We have \[k=ij=i(1j+0k+0)=i(1,0,0)=ie_x.\] Consequently, the spherical components of $k$ are $0$ and $e_x$. An ordered pair $(u, v)$ of unit vectors represents the argument of $k$ if and only if $u \perp e_x$ and $v = e_x \times u$. Given that $e_y \perp e_x$, we select $u = e_y$ and $v = e_x \times e_y = e_z$. Thus, $\alpha_k = (e_y, e_z)$ (\Cref{fig:fig9}). Therefore,
\begin{equation} \label{eqt:eqt6.4}
k = e^{i(e_y, e_z)} \ \text{with} \ \arg(k) = (e_y, e_z).
\end{equation}
\end{enumerate}
\end{example}

\subsection{Geometric interpretation of the formula $\boldsymbol{ki=j}$}\label{subsec:subsec6.2}
In \Cref{exm:exm6.1}, we determined the arguments of the unit quaternions $i$, $j$, and $k$. Utilizing (\ref{eqt:eqt6.2}), (\ref{eqt:eqt6.3}), (\ref{eqt:eqt6.4}), and \Cref{propo:propo5.1}.(1), we obtain the equivalent formulations, \[ki=j \Leftrightarrow e^{i\alpha_k}e^{i\alpha_i}=e^{i\alpha_j} \Leftrightarrow e^{i(\alpha_i+\alpha_k)}=e^{i\alpha_j}.\] Applying (\ref{eqt:eqt5.10}), this implies that $\mu(\alpha_i + \alpha_k) = \mu(\alpha_j)$. Consequently, $\alpha_i + \alpha_k = \alpha_j$, leading to the relation $(e_x, e_y) + (e_y, e_z) = (e_x, e_z)$.

\begin{figure}[H]%
\centering
\includegraphics[width=1\textwidth]{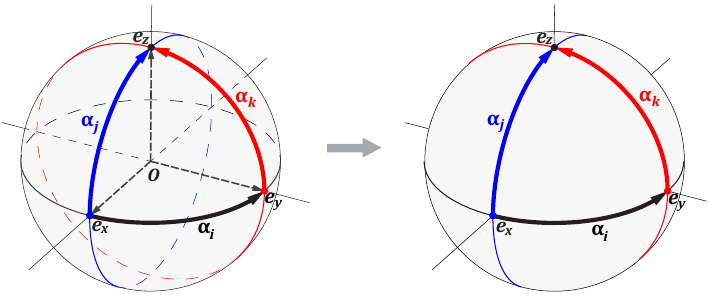}
\caption{Geometric interpretation of the equation $ki = j$ through the Chasles relation applied to the arguments of quaternions $i$, $j$, and $k$.}\label{fig:fig9}
\end{figure}

As illustrated in \Cref{fig:fig9}, the multiplication relation $ki = j$ can be geometrically interpreted through the straightforward Chasles relation, \[(e_x,e_y)+(e_y,e_z)=(e_x,e_z).\]

\subsection{Geometric interpretation of the non-commutativity of unit quaternion multiplication}\label{subsec:subsec6.3}
\begin{figure}[H]%
\centering
\includegraphics[width=1\textwidth]{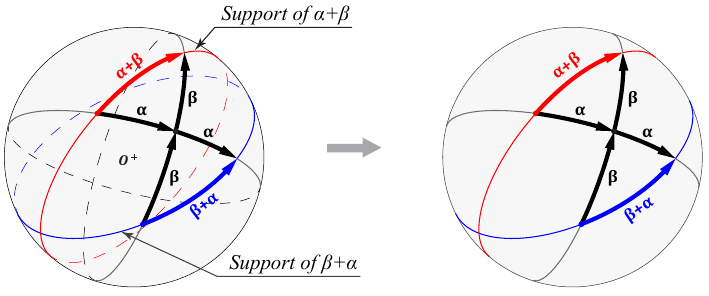}
\caption{The two arguments $\alpha + \beta$ and $\beta + \alpha$ possess identical scalar components; however, they do not share the same support.}
\label{fig:fig10}
\end{figure}

Let $p$ and $q$ be two unit quaternions with respective arguments $\alpha$ and $\beta$, where $\alpha$ and $\beta$ are spherical-vectors. According to \Cref{propo:propo5.2}, we obtain $\arg(pq) = \beta + \alpha$ and $\arg(qp) = \alpha + \beta$. Consequently, the non-commutativity of the multiplication of $p$ and $q$ can be geometrically interpreted through the non-commutativity of the addition of their arguments $\alpha$ and $\beta$, as depicted in \Cref{fig:fig10} and further demonstrated in the following numerical \Cref{exm:exm6.2}.

\begin{example}\label{exm:exm6.2}
Consider the following two unit quaternions $p = \sqrt{6}/6(2-j-k)$ and $q = \sqrt{2}/2(1+i)$. In this final example, we employ spherical-vectors to provide a geometric representation for the quaternions $p$ and $q$, as well as their products $qp = \sqrt{3}/3(1+i-k)$ and $pq = \sqrt{3}/3(1+i-j)$. Let us denote the respective arguments of the quaternions $p$, $q$, $qp$, and $pq$ as $\alpha_{p}$, $\alpha_{q}$, $\alpha_{qp}$, and $\alpha_{pq}$. According to \Cref{propo:propo5.2}, $\alpha_{qp} = \alpha_p + \alpha_q$ and $\alpha_{pq} = \alpha_q + \alpha_p$. To represent these four spherical-vectors on the unit sphere $S^2$, we first need to express $p$ and $q$ in their spherical forms. To do so, we utilize \Cref{def:def3.2} and \Cref{rem:rem3.1}:
\[p=\frac{\sqrt{6}}{3}+i\frac{\sqrt{6}}{6}\left(-j+k+0\right)=\frac{\sqrt{6}}{3}+i\left(-\frac{\sqrt{6}}{6},\frac{\sqrt{6}}{6},0\right)\] and \[q=\frac{\sqrt{2}}{2}+i\left(0j+0k+\frac{\sqrt{2}}{2}\right)=\frac{\sqrt{2}}{2}+i\left(0,0,\frac{\sqrt{2}}{2}\right).\]
The spherical components of $p$ are given by \[\lambda_{p} = \frac{\sqrt{6}}{3} \text{ and } n_{p} = \left(-\frac{\sqrt{6}}{6}, \frac{\sqrt{6}}{6}, 0\right),\] while the spherical components of $q$ are \[\lambda_{q} = \frac{\sqrt{2}}{2} \text{ and } n_{q} = \left(0, 0, \frac{\sqrt{2}}{2}\right).\] According to \Cref{rem:rem5.1}.(i), the components of $\alpha_{p}$ are the spherical components of $p$, and the components of $\alpha_{q}$ are the spherical components of $q$. Then, as per \Cref{lem:lem4.2}.(i), given that $n_{p} \times n_{q}$ is non-zero, there are only two possible combinations of three unit vectors $u$, $v$, and $w$ such that $\alpha_{p} = (u, v)$ and $\alpha_{q} = (v, w)$. We may choose, for instance, the combination \[\begin{cases} v=\frac{n_{p} \times n_{q}}{\Vert n_{p} \times n_{q}\Vert}\\ u=\lambda_{p}v-n_{p}\times v\\ w=\lambda_{q}v+n_{q}\times v. \end{cases}\]
Upon calculation, we obtain the following results: \[\begin{cases} v=\left(\frac{\sqrt{2}}{2},\frac{\sqrt{2}}{2},0\right)\\ u=\left(\frac{\sqrt{3}}{3},\frac{\sqrt{3}}{3},\frac{\sqrt{3}}{3}\right)\\ w=\left(0,1,0\right)=e_y. \end{cases}\] Consequently, we can represent the product $qp$ through the sum $\alpha_{p} + \alpha_{q} = \alpha_{qp}$, as depicted in \Cref{fig:fig11}.(1).\\

Similarly, we observe that $n_{q} \times n_{p} \neq \vec{0}$, indicating that there are two combinations of three unit vectors $u'$, $v'$, and $w'$ such that $\alpha_q = (u', v')$ and $\alpha_p = (v', w')$. We choose, for instance, the combination:
\[\begin{cases} v'=\frac{n_{p} \times n_{q}}{\Vert n_{p} \times n_{q}\Vert}=v\\ u'=\lambda_{q}v-n_{q}\times v\\ w'=\lambda_{p}v+n_{p}\times v. \end{cases}\] Which results in: \[\begin{cases} v'=v=\left(\frac{\sqrt{2}}{2},\frac{\sqrt{2}}{2},0\right)\\ u'=\left(1,0,0\right)=e_x\\ w'=\left(\frac{\sqrt{3}}{3},\frac{\sqrt{3}}{3},-\frac{\sqrt{3}}{3}\right).\end{cases}\] We can now represent the product $pq$ through the sum $\alpha_{q} + \alpha_{p} = \alpha_{pq}$, as depicted in \Cref{fig:fig11}.(1).

\begin{figure}[H]%
\centering
\includegraphics[width=1\textwidth]{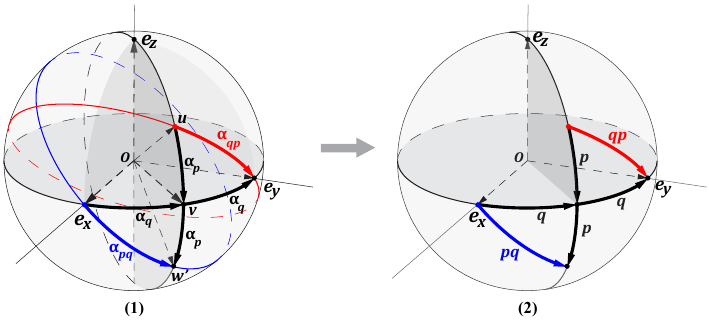}
\caption{$-$ Illustration 1: Geometric interpretation of the products $qp$ and $pq$ via the Chasles relations $\alpha_p + \alpha_q = \alpha_{qp}$ and $\alpha_q + \alpha_p = \alpha_{pq}$, respectively.\\
$-$ Illustration 2: Establishing correspondence between unit quaternions $p$, $q$, $qp$, $pq$ and their respective arguments $\alpha_{p}$, $\alpha_{q}$, $\alpha_{qp}$, and $\alpha_{pq}$.}
\label{fig:fig11}
\end{figure}

Due to the identification established throughout this paper between unit quaternions and spherical-vectors, we now possess a robust geometric approach for representing unit quaternions (which belong to a 4-dimensional space) and their products on a surface of a 3-dimensional space. As observed in (\Cref{fig:fig11}.(2)), we adopt the notation $p = \alpha_p$, $q = \alpha_q$, $qp = \alpha_p + \alpha_q$, and $pq = \alpha_q + \alpha_p$.
\end{example}

\textbf{Acknowledgement}. We extend our gratitude to the \textit{les-mathematiques.net} forum, its members, and particularly to the anonymous professor GaBuZoMeu, whose assistance was instrumental in the formalization of \Cref{def:def2.2} pertaining to spherical-vectors.

\newpage

%

\bibliographystyle{amsplain}  

\providecommand{\bysame}{\leavevmode\hbox to3em{\hrulefill}\thinspace}

\end{document}